\documentclass[10pt,conference]{IEEEtran}

\usepackage[english]{babel}
\usepackage{graphicx}
\usepackage{subcaption}
\usepackage{framed}
\usepackage[normalem]{ulem}
\usepackage{amsmath}
\usepackage{amsthm}
\usepackage{amssymb}
\usepackage{amsfonts}
\usepackage{enumerate}
\usepackage{algorithmic}
\usepackage[utf8]{inputenc}
\usepackage[top=1 in,bottom=1in, left=1 in, right=1 in]{geometry}
\usepackage{wrapfig}
\usepackage[ruled,vlined]{algorithm2e}
\usepackage{tikz}
\usepackage{float}

\usetikzlibrary{positioning}
\newtheorem{question}{Question}
\newtheorem{theorem}{Theorem}
\newtheorem{definition}{Definition}
\newtheorem{example}{Example}

\newcommand{\RR}{\mbox{$\mathbb{R}$}} 
\newcommand{\tens}[1]{\mathcal{#1}}

\newcommand{\matr}[1]{\mathbf{#1}}

\newcommand{\kron}{\otimes}  
\newcommand{\khatri}{\odot} 
\newcommand{\hadam}{\boxdot} 
\newcommand{\out}{\pmb\otimes}  
\makeatletter
\def\vecd{\mathop{\operator@font vecd}\nolimits}
\newcommand{\trace}[1]{\mathop{\operator@font trace}\{#1\}}
\newcommand{\Diag}[1]{\mathop{\operator@font Diag}\{#1\}} 
\newcommand{\diag}[1]{\mathop{\operator@font diag}\{#1\}} 

\newcommand{\rank}[1]{\mathop{\operator@font rank}\{#1\}}
\def\vecd{\mathop{\operator@font vecd}\nolimits}
\newcommand{\krank}[1]{\mathop{\operator@font krank}\{#1\}}

\makeatletter
\newcommand{\removelatexerror}{\let\@latex@error\@gobble}
\makeatother
\usepackage{multirow}

\title{Tensor Completion through Total Variation with Initialization from Weighted HOSVD}
\author{Zehan Chao, Longxiu Huang and Deanna Needell}
\date{}

\begin{document}

\maketitle

\begin{abstract}
In our paper, we have studied the tensor completion problem when the sampling pattern is deterministic. We first propose a simple but efficient weighted HOSVD algorithm for recovery from noisy observations. Then we use the weighted HOSVD result as an initialization for the total variation. We have proved the accuracy of the weighted HOSVD algorithm from theoretical and numerical perspectives. In the numerical simulation parts, we also showed that by using the proposed initialization, the total variation algorithm can efficiently fill the missing data for images and videos.
\end{abstract}
\section{Introduction}
Tensor, a high-dimensional array which is an extension of matrix,  plays an important role in a wide range of real world applications \cite{AD1981,comon2000}. Due to the high-dimensional structure, tensor could preserve more information compared to the unfolded matrix. For instance, a $k$ frame, $m\times n$ video stored as an $m\times n \times k$ tensor will keep the connection between each frames, splitting the frames or unfolding this tensor may lose some conjunctional information. A brain MRI would also benefit from the 3D structure if stored as tensor instead of 
randomly arranging several snapshots as matrices. On the other hand, most of the real world datasets are partially missing and incomplete data which can lead to extremely low performance of downstream applications. The linear dependency and redundancy between missing and existing data can be leveraged to recover unavailable data and improve the quality and scale of the incomplete dataset. The task of recovering missing elements from partially observed tensor is called tensor completion and has attracted widespread attention in many applications.
e.g., image/video inpainting \cite{kressner2014low,liu2012tensor}, recommendation systems \cite{symeonidis2008tag}. Matrix completion problem, \cite{ABEV6,AFSU7,CCS10,CZ2016, CP2010,CR9,CS2004,foucart2019weighted,GL2011,GNOT1992} as a special case of tensor completion problem has been well-studied in the past few decades, which enlightened researchers on developing further tensor completion algorithms. Among different types of data matrices, image data is commonly studied and widely used for performance indicator. One traditional way to target image denoising problem is to minimize the total variation norm. Such method is based on the assumption of locally smoothness pattern of the data. Yet in recent decades, thanks to the algorithm development of non-negative matrix factorization (NMF) and nuclear norm minimization (NNM), the low-rank structure assumption becomes increasingly popular and extensively applied in related studies. In both matrix completion and tensor completion studies, researchers are trying to utilize and balance both assumption in order to improve the performance of image recovery and video recovery tasks.

 In this research, we will provide an improved version of total variation minimization problem by providing a proper initialization. To implement the initialization, we have designed a simple but efficient algorithm which we call weighted HOSVD algorithm for low-rank tensor completion from  a deterministic sampling
 pattern, which is motivated from \cite{CJN2019, foucart2019weighted}.

\section{Tensor Completion Problem}
In this section, we would provide a formal definition for the tensor completion problem. First of all, we will introduce notations, basic operations and definitions for tensor.
\subsection{Preliminaries and Notations}
~

\noindent
 Tensors, matrices, vectors and scalars are denoted in different typeface for clarity below. Throughout this paper, calligraphic boldface capital letters are used for tensors,  capital letters are used for matrices, lower boldface letters for vectors, and regular letters for scalars.  The set of the first $d$ natural numbers will be denoted by $[d]:=\{1,\cdots,d\}$. Let $\mathcal{T}\in\mathbb{R}^{d_1\times\cdots\times d_n}$ and $\alpha\in\mathbb{R}$, $\mathcal{T}^{(\alpha)}$ represents the pointwise power operator i.e., $(\mathcal{T}^{(\alpha)})_{i_1\cdots i_n}=(\mathcal{T}_{i_1\cdots i_n})^{\alpha}$. We use $\mathcal{T}\succ 0$ to denote the tensor with $\mathcal{T}_{i_1\cdots i_n}>0$ for all $i_1,\cdots, i_n$.  $\boldsymbol{1}_\Omega$ denotes the tensor with all entries equal to $1$ on $\Omega$ and $0$ otherwise.
\begin{definition}[Tensor]
A tensor is a multidimensional array. The dimension of a tensor is called the order (also called the mode). The space of a real tensor of order $n$ and of size $d_1\times\cdots\times d_n$ is denoted as $\mathbb{R}^{d_1\times  \cdots\times d_n}$. The elements of a tensor $\mathcal{T}\in\mathbb{R}^{d_1\times \cdots\times d_n}$ are denoted by $\mathcal{T}_{i_1\cdots i_n}$.
\end{definition}
For an order $n$ tensor $\mathcal{T}$ can be matricized in $n$ ways by unfolding it along each of the $n$ modes, next we will give the definition for the matricization of a given tensor.
\begin{definition}[Matricization of a tensor] 
The mode-$k$ matricization of tensor $\mathcal{T}\in\mathbb{R}^{d_1\times \cdots\times d_n}$ is the matrix, which is denoted as $\mathcal{T}_{(k)}\in\mathbb{R}^{d_k\times\prod_{j\neq k}d_j}$, whose columns are composed of all the vectors obtained from $\mathcal{T}$ by fixing all indices but $i$th.
\end{definition}
In order to illustrate the matricization of a tensor, let us consider the following example.
\begin{example}
Let $\mathcal{T}\in\mathbb{R}^{3\times 4\times 2}$ with the following frontal slices:
\[T_1=\begin{bmatrix}
1&4&7&10\\
2&5&8&11\\
3&6&9&12
\end{bmatrix} \quad~~~ T_2=\begin{bmatrix}
13&16&19&22\\
14&17&20&23\\
15&18&21&24
\end{bmatrix},
\]
then the three mode-n matricizations  are
\begin{eqnarray*}
\mathcal{T}_{(1)}&=&\begin{bmatrix}
1&4&7&10&13&16&19&22\\
2&5&8&11&14&17&20&23\\
3&6&9&12&15&18&21&24
\end{bmatrix},\\
\mathcal{T}_{(2)}&=&\begin{bmatrix}
1&2&3&13&14&15\\
4&5&6&16&17&18\\
7&8&9&19&20&21\\
10&11&12&22&23&24
\end{bmatrix},\\
\mathcal{T}_{(3)}&=&\begin{bmatrix}
1&2&3&\cdots&10&11&12\\
13&14&15&\cdots&22&23&24
\end{bmatrix}.
\end{eqnarray*}
\end{example}
\begin{definition}[Folding Operator]
Suppose $\mathcal{T}$ be a tensor. The mode-$k$ folding operator of a matrix $M=\mathcal{T}_{(k)}$, denoted as $\text{fold}_k(M)$, is the inverse operator of the unfolding operator.
\end{definition}

\begin{definition}[$\infty$ norm]
Let $\mathcal{T}\in\mathbb{R}^{d_1\times d_2\times\cdots\times d_n}$, the $\|\mathcal{T}\|_\infty$ is defined as
\[\|\mathcal{T}\|_{\infty}=\max_{i_1,i_2,\cdots,,i_n}|\mathcal{T}_{i_1i_2\cdots i_n}|.
\] The unit ball under $\infty$ norm is denoted by $\boldsymbol{B}_\infty$.
\end{definition}
\begin{definition}[Frobenius norm]
The Frobenius norm for tensor $\mathcal{T}\in\mathbb{R}^{d_1\times d_2\times \cdots\times d_n}$ is defined as
\[\|\mathcal{T}\|_F=\sqrt{\sum_{i_1,i_2,\cdots, i_n}\mathcal{T}_{i_1i_2\cdots i_n}^2}.\]
\end{definition}
\begin{definition}[Product Operations]~
\begin{enumerate}
\item[$\bullet$] Mode-$k$ Product: Mode-$k$ product of tensor  $\mathcal{T}\in\mathbb{R}^{d_1\times\cdots\times d_n}$ and matrix $A\in\mathbb{R}^{d\times d_k}$ is defined by
\[\mathcal{T}\times_k A=\text{fold}_k(A\mathcal{T}_{(k)}),
\]
i.e.,
\[(\mathcal{T}\times_k A)_{i_1\cdots i_{k-1}j i_{k+1}\cdots i_n}=\sum_{i_k=1}^{d_k}\mathcal{T}_{i_1i_2\cdots i_n}A_{ji_k}.
\]
\item[$\bullet$]Outer product: Let $\boldsymbol{a_1}\in\mathbb{R}^{d_1}, \cdots, \boldsymbol{a_n}\in\mathbb{R}^{d_n}$. The outer product among these $n$ vectors is a tensor $\mathcal{T}\in\mathbb{R}^{d_1\times\cdots\times d_n}$ defined as:
\[\mathcal{T}=\boldsymbol{a_1}\out \cdots\out \boldsymbol{a_n}, ~~ \mathcal{T}_{i_1,\cdots,i_n}=\prod_{k=1}^n \boldsymbol{a_k}(i_k).\]
    \item [$\bullet$]Kronecker product of matrices: The Kronecker product of $A\in\mathbb{R}^{I\times J}$ and $B\in\mathbb{R}^{K\times L}$ is denoted by $A\kron B$. The result is a matrix of size $(KI)\times (JL)$ and defined by
\begin{eqnarray*}
A\kron B&=&\begin{bmatrix}
A_{11}B&A_{12}B&\cdots &A_{1J}B\\
A_{21}B&A_{22}B&\cdots &A_{2J}B\\
\vdots&\vdots&\ddots&\vdots\\
A_{I1}B&A_{I2}B&\cdots&A_{IJ}B
\end{bmatrix}.
\end{eqnarray*}
\item[$\bullet$] Khatri-Rao product:
Given matrices $A\in\mathbb{R}^{d_1\times r}$ and $B\in\mathbb{R}^{d_2\times r}$, their Khatri-Rao product is denoted by $A\khatri B$. The result is a matrix of size  $(d_1d_2)\times r$ defined by
\[A\khatri B=\begin{bmatrix}
\boldsymbol{a_1}\kron \boldsymbol{b_1}&\cdots&\boldsymbol{a_r}\kron\boldsymbol{ b_r}
\end{bmatrix},
\]
where $\boldsymbol{a_i}$ and $\boldsymbol{b_i}$ stands for the $i$-th column of $A$ and $B$ respectively. 
\item [$\bullet$]Hadamard product: Given two tensors $\mathcal{T}$ and $\mathcal{Y}$, both of size $d_1\times \cdots\times d_n$, their Hadamard product is denoted by $\mathcal{X}\hadam \mathcal{Y}$. The result is also of the size $d_1\times d_2\times \cdots \times d_n$ and the elements of $\mathcal{X}\hadam\mathcal{Y}$ are defined as the elementwise tensor product i.e.,
\[(\mathcal{X}\hadam\mathcal{Y})_{i_1i_2\cdots i_n}=\mathcal{X}_{i_1i_2\cdots i_n}\mathcal{Y}_{i_1i_2\cdots i_n}.
\]

\end{enumerate}
\end{definition}
\begin{definition}[Rank-one Tensors]
An $n$-order tensor $\mathcal{T}\in\mathbb{R}^{d_1\times d_2\times\cdots\times d_n}$ is rank one if it can be written as the out product of $n$ vectors, i.e.,
\[\mathcal{T}=\boldsymbol{a_1}\out\cdots\out \boldsymbol{a_n}.\] 
\end{definition}
\begin{definition}[Tensor (CP) rank\cite{HF1927,HF1928}]
The \textit{rank} of a tensor $\mathcal{X}$, denoted $\text{rank}(\mathcal{X})$, is defined as the smallest number of rank-one tensors that generate $\mathcal{X}$ as their sum. We use $K_r$ to denote the cone of rank-r tensors.
\end{definition}

Different from the case of matrices, the rank of a tensor is presently not understood well. And the problem of computing the rank of a tensor is NP-hard.  Next we will introduce a new rank definition related to the tensor.
\begin{definition}[Tucker rank  \cite{HF1928}]Let $\mathcal{X}\in\mathbb{R}^{d_1\times\cdots\times d_n}$.
The tuple $(r_1,\cdots,r_n)\in\mathbb{N}^{n}$, where $r_k=\text{rank}(\mathcal{X}_{(k)})$ is called tensor Tucker rank of $\mathcal{X}$. We use $K_{\boldsymbol{r}}$ to denote the cone of Tucker rank $\boldsymbol{r}$ tensors.
\end{definition}

\subsection{Problem Statement}
In order to find a good initialization for total variation method, we would like to solve the following questions.
\begin{question}
Given a deterministic sampling pattern $\Omega$ and corresponding (possibly noisy) observations from the tensor, what type of recovery error can we expect, in what metric, and how may we efficiently implement this recovery?
\end{question}
\begin{question}
Given a sampling pattern $\Omega$, and noisy observations $\mathcal{T}_\Omega+\mathcal{Z}_\Omega$, for what rank-one weight tensor $\mathcal{H}$ can we efficiently find a tensor $\widehat{\mathcal{T}}$ so that $\|\mathcal{H}\hadam(\widehat{\mathcal{T}}-\mathcal{T})\|_F$ is small compared to $\|\mathcal{H}\|_F$? And how can we efficiently find such weight tensor $\mathcal{H}$, or certify that a fixed $\mathcal{H}$ has this property?
\end{question}
In order to find weight tensor, we consider the optimization problem
\begin{eqnarray*}
\mathcal{W}:=\underset{{\mathcal{X}\succ 0, \text{rank}(\mathcal{X})=1}}{\text{argmin}}\|\mathcal{X}-\boldsymbol{1}_{\Omega}\|_F
\end{eqnarray*}
$\mathcal{W}$ can be estimated by using the least square CP algorithm \cite{CC1970,HR1970}. After we find $\mathcal{W}$, then we consider the following optimization problem to estimate $\mathcal{T}$:
\begin{equation}\label{eqn:obj1}
    \widehat{\mathcal{T}}=\mathcal{W}^{(-1/2)}\hadam\underset{\text{Tucker\_rank}(\mathcal{T})=\boldsymbol{r}}{\text{argmin}}\|\mathcal{T}-\mathcal{W}^{(-1/2)}\hadam\mathcal{Y}_\Omega \|_F,
\end{equation}
where $\mathcal{Y}_\Omega=\mathcal{T}_\Omega+\mathcal{Z}_\Omega$. As we know, to solve problem \eqref{eqn:obj1} is NP-hard. 
In order to solve \eqref{eqn:obj1} in polynomial time, we consider the HOSVD process \cite{DDV2000}. Assume that $\mathcal{T}$ has Tucker rank $\boldsymbol{r}=[r_1,\cdots,r_n]$.  Let
\[\widehat{A}_i=\underset{\text{rank}(A)=r_i}{\text{argmin}}\|A-(\mathcal{W}^{(-1/2)}\hadam \mathcal{Y}_\Omega)_{(i)}\|_2.
\]
and set $\widehat{U}_i$ to be the left singular vector matrix of $\widehat{A}_i$. Then the estimated tensor is of the form \[\widehat{\mathcal{T}}=\mathcal{W}^{(-1/2)}\hadam((\mathcal{W}^{(-1/2)}\hadam\mathcal{Y}_\Omega)\times_1 \widehat{U}_1\widehat{U}_1^T\times_2 \cdots\times_n \widehat{U}_n\widehat{U}_n^T.\]
 In the following, we call our algorithm weighted HOSVD algorithm.

With the output from weighted HOSVD, $\mathcal{\widehat{T}}$, we can solve the following total variation problem with $\mathcal{\widehat{T}}$ as initialization:

\begin{align*}
    \min_\tens{X}\quad & \|\tens{X}\|_{TV} 
    \\
    s.t. \quad &\mathcal{X}_\Omega=\mathcal{Y}_\Omega.
\end{align*}

This total variation minimization problem is solved by iterative method. We will discuss the details of algorithm and numerical performance in section IV and V.

\section{Related Work}
In this section, we will briefly step through the history of matrix completion and introduce several relevant studies on tensor completion.

Suppose we have a partially observed matrix $M$ under the low-rank assumption and give a deterministic sampling pattern $\Omega$, the most intuitive optimization problem raised here is:
\vspace*{-0.1cm}
\begin{align*}
    \min_X\quad & Rank(X), \\
    s.t. \quad &  X_\Omega=M_\Omega.
\end{align*}

However, due to the computational complexity (NP-hard) of this minimization problem, researchers developed workarounds by defining new optimization problems which could be done in polynomial time. Two prominent substitutions are the nuclear norm minimization(NNM) and low rank matrix factorization(LRMF):
\vspace*{-0.1cm}
\begin{align*}
    \min_X\quad & \|X\|_* \\
    s.t. \quad &  X_\Omega=M_\Omega,
\end{align*}
where $\|\cdot\|_*$ stands for the sum of singular values. 
\vspace*{-0.1cm}
\begin{align*}
    \min_{A,B}\quad & \|(X-AB)_\Omega\|_F \\
    s.t. \quad &  A\in{\RR^{m\times r}},B\in{\RR^{r\times n}},
\end{align*}
where $A\in{\RR^{m\times r}},B\in{\RR^{r\times n}}$ are restricted to be the low rank components.

While the task went up to tensor completion, the low rank assumption became even harder to approach. Some of the recent studies performed matrix completion algorithms on the unfolded tensor and obtained considerable results. For example, \cite{liu2012tensor} introduced nuclear norm to unfolded tensors and took the weighted average for loss function. They proposed several algorithms such as FaLRTC and SiLRTC to solve the minimization problem:
\begin{align*}
    \min_\tens{X}\quad & \sum_{i=1}^n\alpha_i\|\tens{X}_{(i)}\|_* \\
    s.t. \quad &  \tens{X}_\Omega=\tens{T}_\Omega.
\end{align*}

\cite{xu2013parallel} applied low-rank matrix factorization(LRMF) to  all-mode unfolded tensors and defined the minimization problem as following:
\begin{align*}
    \min_{\tens{X},\matr{A,B}}\quad & \sum_{i=1}^n\alpha_i\|\tens{X}_{(i)}-A_iB_i\|_F^2 \\
    s.t. \quad &  \tens{X}_\Omega=\tens{T}_\Omega
\end{align*}
Where $\matr{A} = \{A_1,...,A_n\}$, $\matr{B} = \{B_1,...,B_n\} $ are the set of low rank matrices with different size according to the unfolded tensor. This method is called TMac and could be solved using alternating minimization.

While researchers often test the performance of their tensor completion algorithms on image/video/MRI data, they started to combine NNM and LRMF with total variation norm minimization when dealing with relevant recovery tasks. For example, \cite{ji2016tensor} introduced the TV regularization into the tensor completion problem:
\vspace*{-0.1cm}
\begin{align*}
    \min_{\tens{X},A,B}\quad & \sum_{i=1}^n\alpha_i\|\tens{X}_{(i)}-A_i B_i\|_F^2+ \mu\|B_3\|_{TV} \\
    s.t. \quad &  \tens{X}_\Omega=\tens{T}_\Omega.
\end{align*}
Note that the $B_3$ here only compute the TV-norm for the first 2 modes. For example, assume that $\tens{X}$ is a video which can be treated as a 3-order tensor, then this TV-norm only counts the variation within each frame without  the variation between frames.



For specific tensors - RGB image data, \cite{li2018total} unfolded the tensor in 2 ways (the horizontal and vertical dimension) and minimized the TV and nuclear norms of each unfolded matrix:
\vspace*{-0.1cm}
\begin{align*}
    \min_\tens{X}\quad & \sum_{i=1}^2(\alpha_i\|\tens{X}_{(i)}\|_* + \mu\|\tens{X}_{(i)}\|_{TV}), \\
    s.t. \quad &  \tens{X}_\Omega=\tens{T}_\Omega.
\end{align*}

In our experiment, we noticed that for a small percentage of observations, for instance, 50\% or more entries are missing, the TV-minimization recovery will produce a similar structure of the original tensor, in the sense of singular values, but NNM will force a large portion of smaller singular values to be zero, which cannot be ignored in the original tensor. Therefore one should be really careful with the choice of minimization problem when performing the completion tasks on a specific dataset. We will discuss the details in the experiment section.

\section{TV Minimization Algorithm}
\subsection{Matrix Denoising Algorithm}
Total Variation Norm is often discretized by \cite{getreuer2012rudin}
\begin{equation*}
\|u\|_{TV}\approx\sum_{i,j}\sqrt{(\nabla_x u)^2_{i,j}+(\nabla_y u)^2_{i,j}}.
\end{equation*}
Hence the image denoising problem is defined as:
\begin{equation*}
\min_{X} \sum_{i,j} \sqrt{(\nabla_x M_{i,j})^2 + (\nabla_y M_{i,j})^2} + \lambda\|M-X\|_F.
\end{equation*}
This can be solved by implementing  the algorithm in \cite{chambolle2010introduction}.


\subsection{Tensor Completion with TV}

Similar to the image denoising algorithm, we first compute the divergence at each entry and move each entry towards the divergence direction. To keep the existing entries unchanged, we project the observed entries to their original values at each step.
 We consider the following minimization problem:
\begin{align*}
\min_{\tens{X}} \quad & \|\tens{X}\|_{TV},\\
s.t. \quad & \tens{X}_\Omega = \tens{T}_\Omega.
\end{align*}
The related algorithm is summarized in Algorithm \ref{alg:TCTV}. 

\removelatexerror
\begin{algorithm*}[H]
\label{alg:TCTV}
\SetKwInOut{Input}{Input}\SetKwInOut{Output}{Output}
\Input{Incomplete tensor $\tens{T}\in\RR^{d_1\times \dots \times d_n}$; Sampling pattern $\Omega\in\{0,1\}^{d_1\times \dots \times d_n}$;
stepsize $h_k$, threshold $\lambda$;
$\tens{X}^0\in\mathbb{R}^{d_1\times\cdots\times d_n}$.}

{Set $\tens{X}^0 = \mathcal{X}^0+(\tens{T}_\Omega-\tens{X}^0_{\Omega})$.}

\For{$k = 0:K$}{
    \For{$i = 1:n$}{
        $\nabla_i(\tens{X}^k_{\alpha_1,...,\alpha_n})=\tens{X}^k_{\alpha_1,...,\alpha_i+1,...,\alpha_n}-\tens{X}^k_{\alpha_1,...\alpha_i,...,\alpha_n},  (\alpha_i=1, 2,...,d_i-1)$
        ($\nabla_i(\cdot) = 0 \text{ when } \alpha_i = d_i$)
        
        $\Delta_i(\tens{X}^k_{\alpha_1,...,\alpha_n})=\tens{X}^k_{\alpha_1,...,\alpha_i-1,...,\alpha_n}+\tens{X}^k_{\alpha_1,...,\alpha_i+1,...,\alpha_n}-2\tens{X}^k_{\alpha_1,...\alpha_i,...,\alpha_n}, (\alpha_i=2,3,...,d_i-1)$
        ($\Delta_i(\cdot) = 0 \text{ when } \alpha_i = 1$ or $d_i$)
    }
    $\Delta(\tens{X}^k_{\alpha_1,...,\alpha_n}) = \sum_i\Delta_i(\tens{X}^k_{\alpha_1,...,\alpha_n})$\\
    $\tens{X}^{k+1}_{\alpha_1,...,\alpha_n} = \tens{X}^k_{\alpha_1,...,\alpha_n}+h_k\cdot\text{shrink}(\frac{\Delta(\tens{X}^k_{\alpha_1,...,\alpha_n})}{\sqrt{\sum_i\nabla_i^2(\tens{X}^k_{\alpha_1,...,\alpha_n})}},\lambda)$\\
    $\tens{X}_\Omega^{k+1} = \tens{T}_\Omega$
    }
\Output{$\tens{X}^K$}
\caption{Tensor Completion through TV minimization}
\end{algorithm*}
\vspace{0.5cm}

In Algorithm 1
, the Laplacian operator computes the divergence of all-dimension gradients for each point. The shrink operator simply moves the input towards 0 with distance $\lambda$, or formally defined as:
\vspace*{-0.03cm}
\[\text{shrink}(x,\lambda) = \mathbf{sign}(x)\cdot \max(|x|-\lambda, 0)\]
\vspace*{-0.03cm}
For $\tens{X}^0$ initialization, simple tensor completion with total variation (TVTC) method would set  $\tens{X}^0$ to be a zero tensor, i.e. $\tens{X}^0 = \mathbf{0}^{d_1\times \dots \times d_n}$, but our proposed method will set $\tens{X}^0$ to be the result from w-HOSVD. We will show the theoretical and experimental advantage of w-HOSVD in the following section.

\section{Main Results}
In order to show the efficiency $\mathcal{\widehat{T}}$ as the initialization for total variation algorithm, we only need to show that $\mathcal{\widehat{T}}$ is close  to $\mathcal{T}$. In the following, the bound of $\|\mathcal{W}\hadam(\mathcal{T}-\mathcal{\widehat{T}})\|_F$ is estimated.
\begin{theorem}\label{thm:gub tensor}
Let $\mathcal{W}=\boldsymbol{w}_1\out\cdots\out\boldsymbol{w}_n \in\mathbb{R}^{d_1\times \cdots\times d_n}$ have strictly positive entries, and fix $\Omega\subseteq[d_1]\times\cdots\times[d_n]$. Suppose that $\mathcal{T}\in\mathbb{R}^{d_1\times \cdots\times d_n}$ has  Tucker-rank $\boldsymbol{r}=[r_1, \cdots,r_n]$ for problem \eqref{eqn:obj1}. Suppose that $\mathcal{Z}_{i_1\cdots i_n}\sim\mathcal{N}(0,\sigma^2)$.
Then with probability at least $1-2^{-|\Omega|/2}$ over the choice of $\mathcal{Z}$,
\begin{multline*}
\|\mathcal{W}^{(1/2)}\hadam(\mathcal{T}-\widehat{\mathcal{T}})\|_F\leq 4\sigma\mu\sqrt{|\Omega|\ln(2)} \\
+2\|\mathcal{T}\|_{\infty}\|\mathcal{W}^{(1/2)}-\mathcal{W}^{(-1/2)}\hadam 1_{\Omega}\|_F,
\end{multline*}
where $\mu^2=\max_{(i_1,\cdots,i_n)\in\Omega}\frac{1}{\mathcal{W}_{i_1\cdots i_n}}$.
\end{theorem}
Notice that the upper bound in Theorem \ref{thm:gub tensor} is for the optimal output $\widehat{\mathcal{T}}$ for problems \eqref{eqn:obj1}, which is general. But the upper bound in Theorem \ref{thm:gub tensor}  contains no    rank information of  the underlying tensor $\mathcal{T}$. In order to introduce the rank information of the underlying tensor $\mathcal{T}$, we restrict our analysis for Problem \eqref{eqn:obj1} by considering the  HOSVD process, we have the following results.
\begin{theorem}[General upper bound for  Tucker rank $\boldsymbol{r}$ tensor]
Let $\mathcal{W}=\boldsymbol{w_1}\out\cdots\out \boldsymbol{w_n}\in\mathbb{R}^{d_1\times\cdots \times d_n}$ have strictly positive entries, and fix $\Omega\subseteq[d_1]\times\cdots\times[d_n]$. Suppose that $\mathcal{T}\in\mathbb{R}^{d_1\times \cdots\times d_n}$ has Tucker rank $\textbf{r}=[r_1~\cdots~ r_n]$. Suppose that $\mathcal{Z}_{i_1\cdots i_n}\sim\mathcal{N}(0,\sigma^2)$ and let 
\[\widehat{\mathcal{T}}=\mathcal{W}^{(-1/2)}\hadam((\mathcal{W}^{(-1/2)}\hadam\mathcal{Y}_\Omega)\times_1 \widehat{U}_1\widehat{U}_1^T\times_2 \cdots\times_n \widehat{U}_n\widehat{U}_n^T)
\]
where $\widehat{U}_1, \cdots, \widehat{U}_n$ are obtained by considering the HOSVD approximation process.
Then with probability at least
\[1-\sum_{i=1}^{n}\frac{1}{d_i+\prod_{j\neq i}d_j}
\]
over the choice of $\mathcal{Z}$, we have
\begin{multline*}\label{uppb}
    \|\mathcal{W}^{(1/2)}\hadam(\mathcal{T}-\widehat{\mathcal{T}})\|_F\lesssim\\ \left(\sum_{k=1}^n \sqrt{r_k\log(d_k+\prod_{j\neq k}d_j)}\mu_k\right)\sigma\\
   +\left(\sum_{k=1}^{n}r_k\|(\mathcal{W}^{(-\frac{1}{2})}\hadam 1_{\Omega}-\mathcal{W}^{(\frac{1}{2})})_{(k)}\|_2 \right)\|\mathcal{T}\|_\infty. 
\end{multline*}
where \begin{multline*}
    \mu_k^2=\max_{i_1,\cdots,i_n}\left\{\sum_{i_1,\cdots,i_{k-1},i_{k+1},\cdots,i_n}\frac{1_{(i_1,i_2,\cdots,i_n)\in\Omega}}{\mathcal{W}_{i_1i_2\cdots i_{n}}},\right.\\
\left.    \sum_{i_k}\frac{1_{(i_1,i_2,\cdots,i_n)\in\Omega}}{\mathcal{W}_{i_1i_2\cdots i_{n}}} \right\}.
    \end{multline*}
\end{theorem}
\section{Simulations}
In this section, we conducted numerical simulations to show the efficiency of the proposed weighted HOSVD algorithm first. Then, we will include the experiment results to show that using the weighted HOSVD algorithm results as a initialization of TV minimization algorithm can accelerate the convergence speed of the original TV minimization. 
\subsection{Simulations for Weighted HOSVD}
In this simulation, we have tested our weighted HOSVD algorithm for 3-order  tensor of the form $\mathcal{T}=\mathcal{C}\times_1 U_1\times_2 U_2\times_3 U_3$ under  uniform and nonuniform   sampling patterns, where $U_i\in\mathbb{R}^{d_i\times r_i}$ and $\mathcal{C}\in\mathbb{R}^{r_1\times r_2\times r_3}$ with $r_i<d_i$. First of all, we generate $\mathcal{T}$ of the size $100\times 100\times 100$ with Tucker rank $\boldsymbol{r}=[r,r,r]$ with $r$ varies from $2$ to $10$. Then we add Gaussian  random noises with $\sigma=10^{-2}$ to $\mathcal{T}$. Next we generate a sampling pattern $\Omega$ which subsample $10\%$ of the entries. We consider estimates $\widehat{\mathcal{T}}_o$, $\widehat{\mathcal{T}}_p$ and $\widehat{\mathcal{T}}_w$ by considering 
\begin{eqnarray*}
\widehat{\mathcal{T}}_o&=&\underset{\text{Tucker\_rank}(\mathcal{X})=\boldsymbol{r}}{\text{argmin}}\|\mathcal{X}-\mathcal{Y}_{\Omega}\|_F,\\
\widehat{\mathcal{T}}_p&=&\underset{\text{Tucker\_rank}(\mathcal{X})=\boldsymbol{r}}{\text{argmin}}\|\mathcal{X}-\frac{1}{p}\mathcal{Y}_{\Omega}\|_F, p=\frac{|\Omega|}{d_1d_2d_3},\\
\widehat{\mathcal{T}}_w&=&\mathcal{W}^{(-1/2)}\hadam\\
&&\underset{\text{Tucker\_rank}(\mathcal{X})=\boldsymbol{r}}{\text{argmin}}\|\mathcal{X}-\mathcal{W}^{(-1/2)}\hadam\mathcal{Y}_{\Omega}\|_F,
\end{eqnarray*}
respectively by using  truncated HOSVD algorithm. 
We give names HOSVD associated with $\widehat{\mathcal{T}}_o$, HOSVD\_p associated with $\widehat{\mathcal{T}}_p$, and w\_HOSVD associated with $\widehat{\mathcal{T}}_w$. For an estimate $\widehat{\mathcal{T}}$, we consider both the weighted and unweighted relative errors:
\[\frac{\|\mathcal{W}^{(1/2)}\hadam(\widehat{\mathcal{T}}-\mathcal{T})\|_F}{\|\mathcal{W}^{(1/2)}\hadam\mathcal{T}\|_F} \text{ and }\frac{\|\widehat{\mathcal{T}}-\mathcal{T}\|_F}{\|\mathcal{T}\|_F}.
\]
We average each error over 20 experiments. 

\begin{figure}
    \centering
    \begin{subfigure}[b]{\linewidth}
		\includegraphics[width=\textwidth]
		{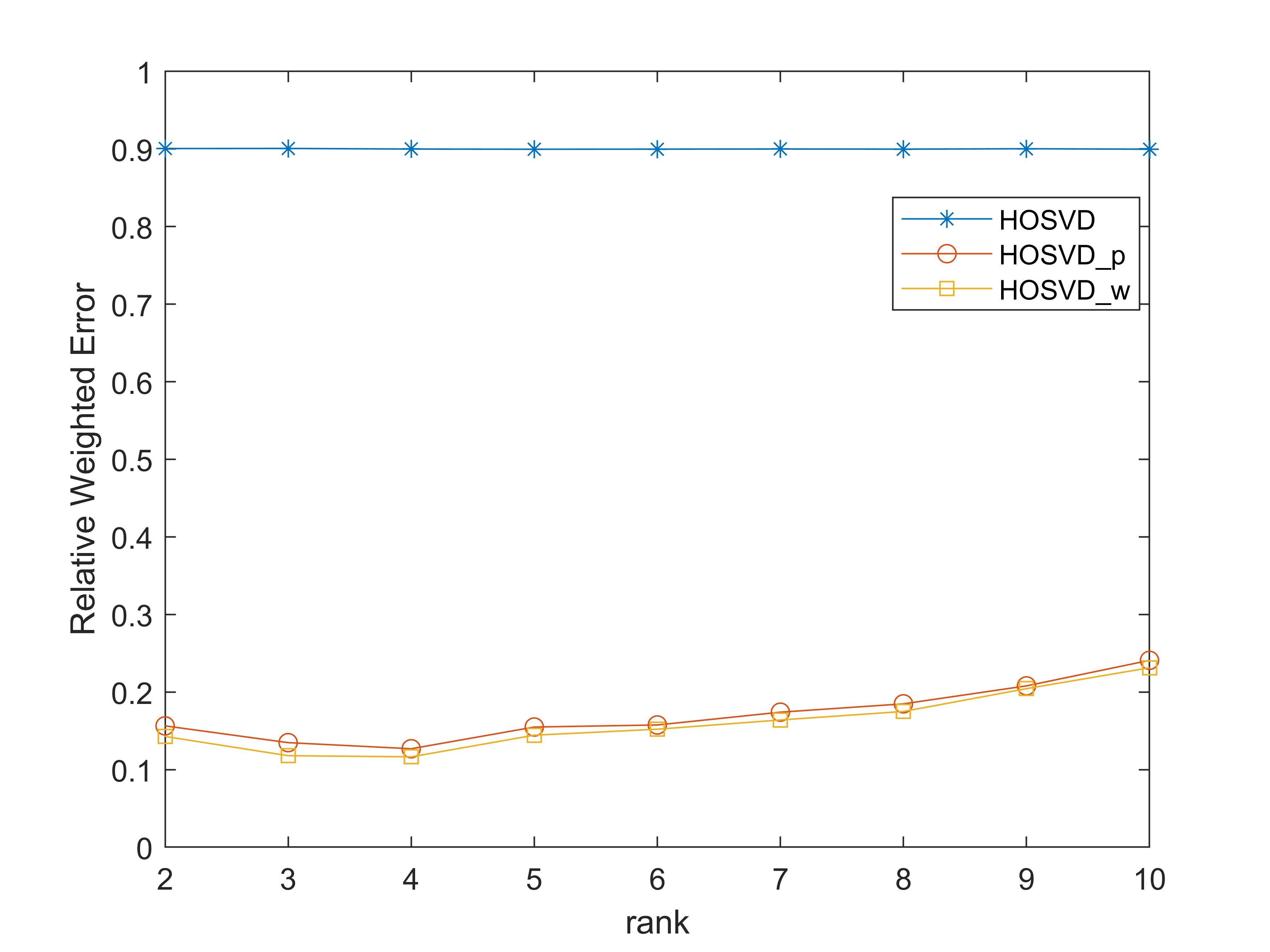}
		\caption{}
	\end{subfigure}
		\begin{subfigure}[b]{\linewidth}
		\includegraphics[width=\textwidth]
		{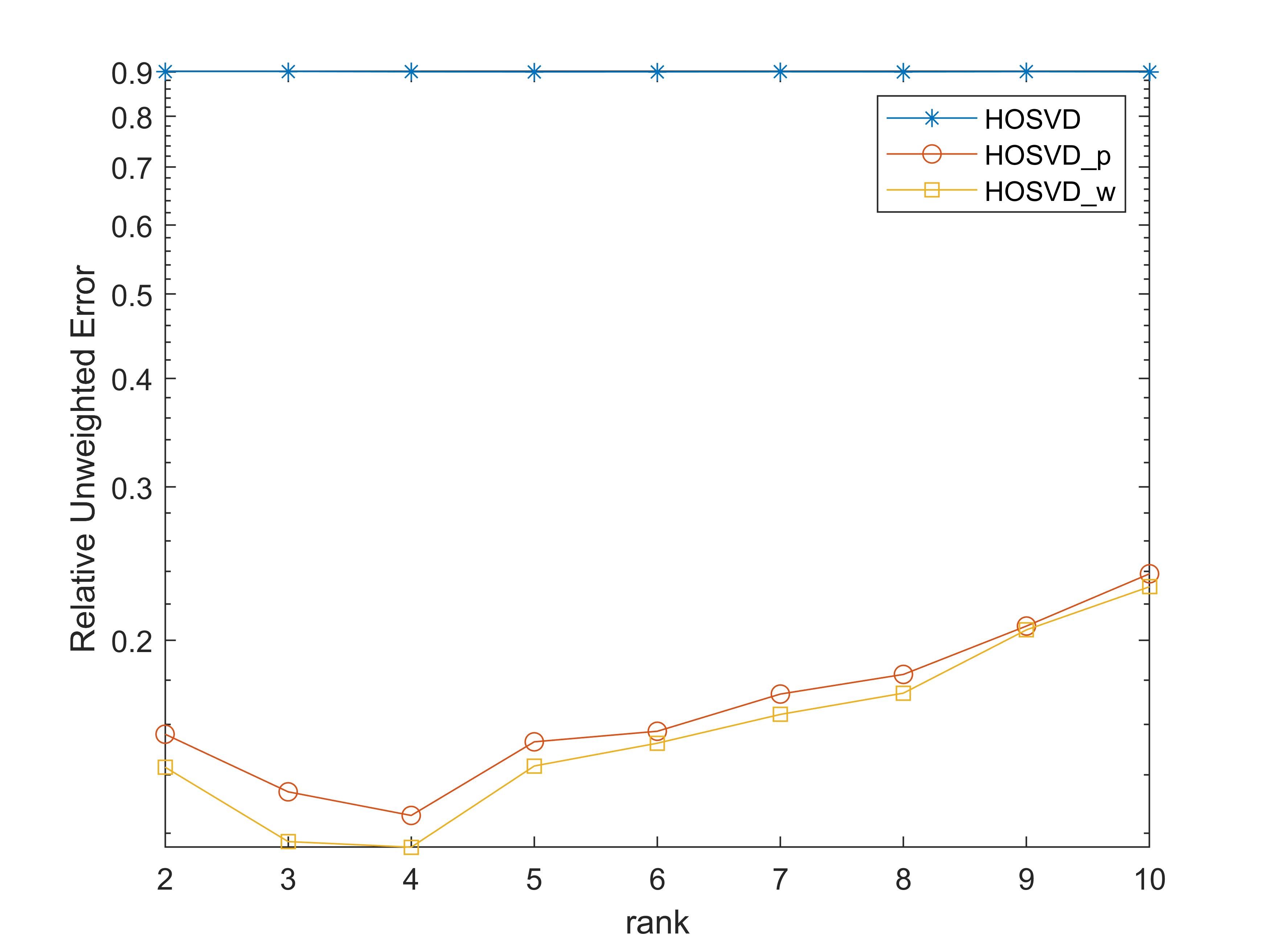}
		\caption{}
	\end{subfigure}
	  \caption{Uniformly subsampled.}
	  \label{FIG: URS}
\end{figure}

When each entry is sampled randomly uniformly with probability $p$, then we have $\mathbb{E}(\mathcal{Y}_\Omega)=p\mathcal{T}$ which implies that the estimate $\widehat{\mathcal{T}}_p$ should perform better than $\widehat{\mathcal{T}}_o$. In Figure \ref{FIG: URS}, we take the sampling pattern to be uniform at random. The estimates $\widehat{\mathcal{T}}_p$ and $\widehat{\mathcal{T}}_w$ perform significantly better than $\widehat{\mathcal{T}}_o$ as expected.

\begin{figure}
    \centering
    	\begin{subfigure}[b]{\linewidth}
		\includegraphics[width=\textwidth]
		{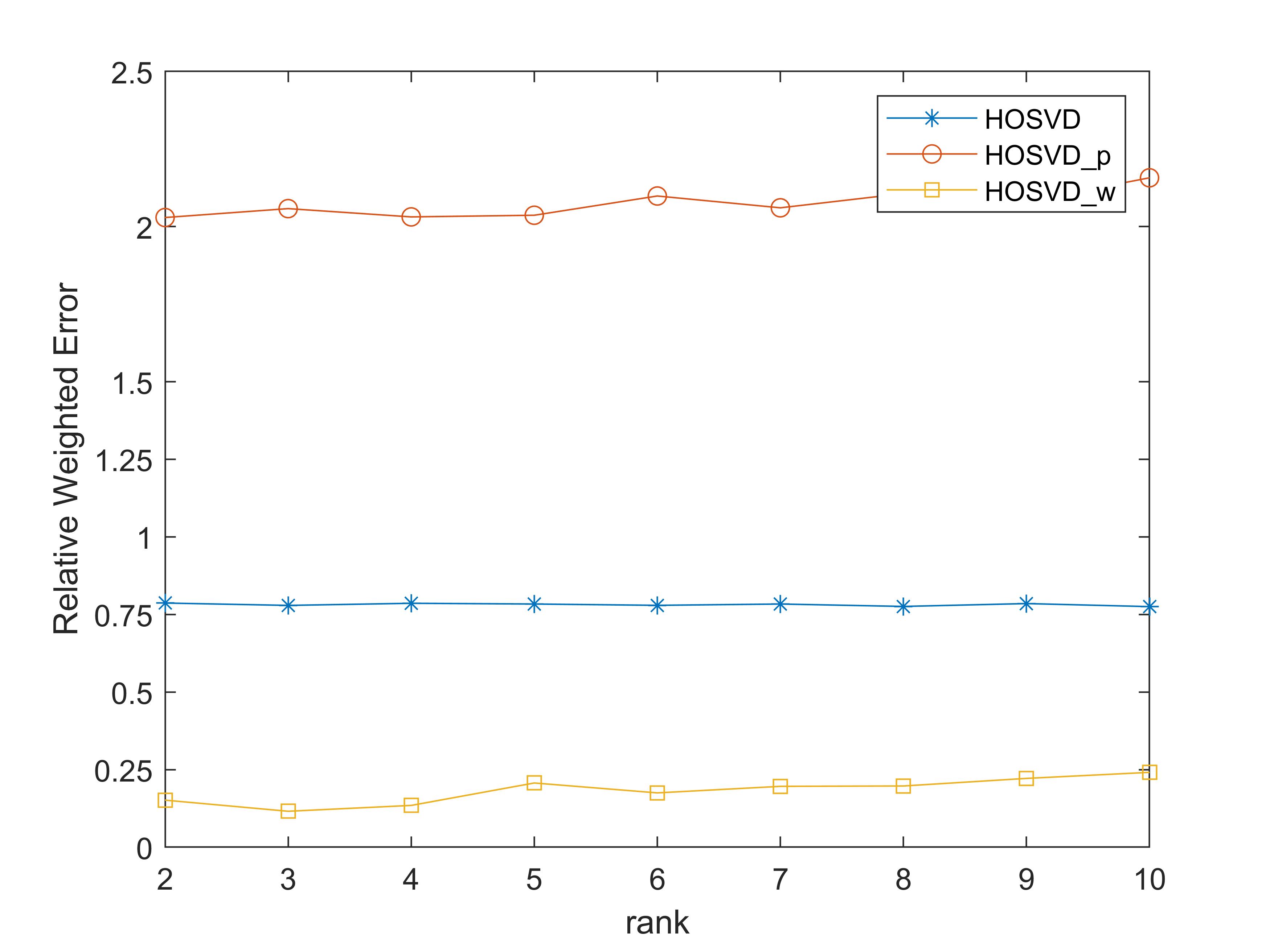}
		\caption{}
	\end{subfigure}
		\begin{subfigure}[b]{\linewidth}
		\includegraphics[width=\textwidth]
		{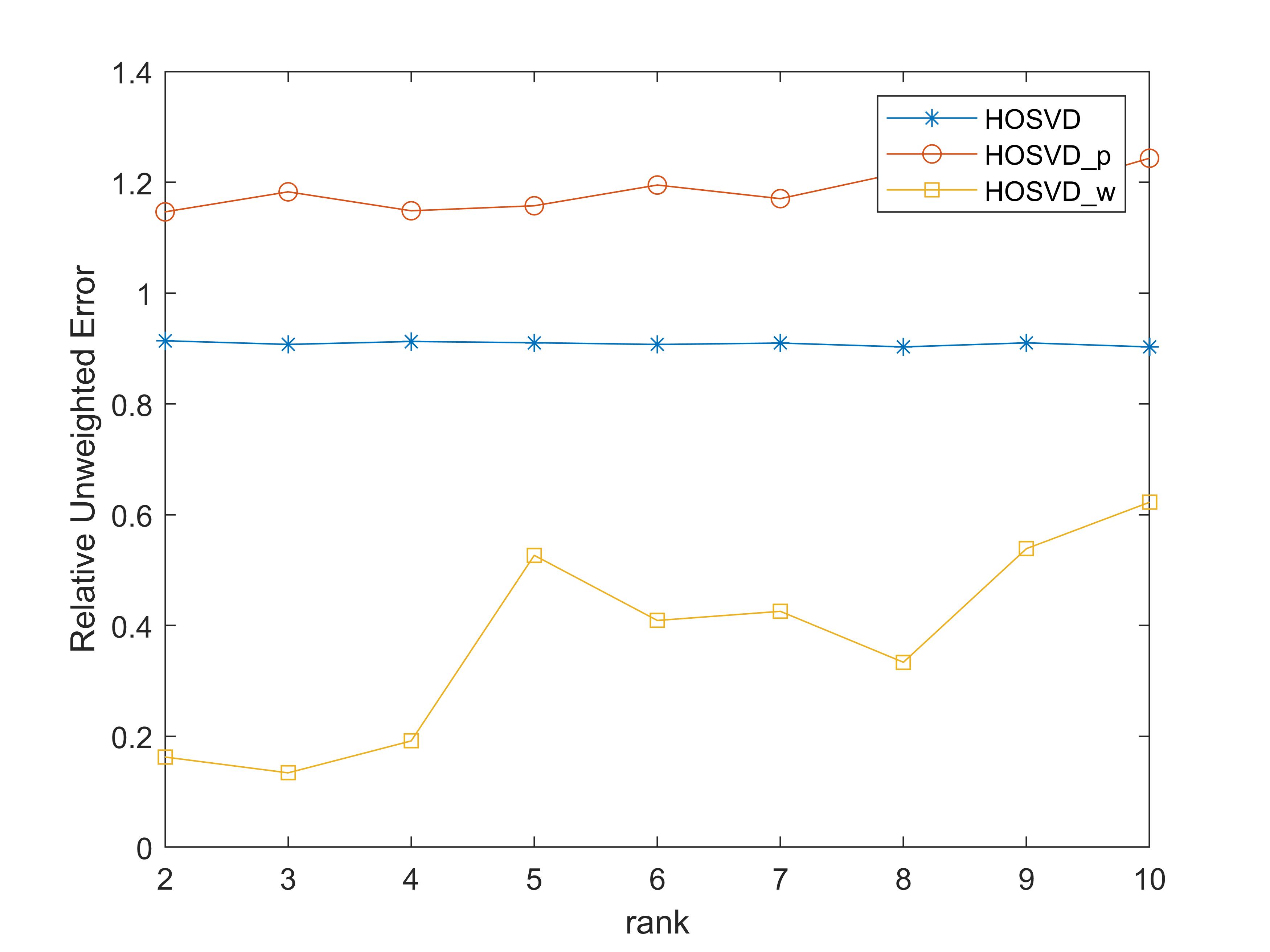}
		\caption{}
	\end{subfigure}
	  \caption{Non-uniformly subsampled.}
	  \label{FIG: NUS}
\end{figure}

    In Figure \ref{FIG: NUS}, the set-up is the same as the one in Figure \ref{FIG: URS} except that the sampling pattern is non-uniformly at random. Then using $\frac{1}{p}$ is not a good weight tensor which as shown in Figure \ref{FIG: NUS}, $\widehat{\mathcal{W}}_p$ works terrible. But $\widehat{\mathcal{T}}_w$ still works better that $\widehat{\mathcal{T}}_o$.

\subsection{Simulations for TV with Initialization from Weighted HOSVD (wHOSVD-TV)}
\subsubsection{Experimental Setup}
In order to show the advantage of weighted HOSVD, we test our proposed algorithm and simple TV minimization method, along with a baseline algorithm on video data. We 
mask a specific ratio of entries and conduct each completion algorithm in order to obtain the completion results $\widehat{\tens{T}}$. The tested sampling rates (SR) are $0.1, 0.3, 0.5 ,0.8$. We then compute the relative root mean square error (RSE):\[
RSE = \frac{\|\widehat{\tens{T}}-\tens{T}\|_F}{\|\tens{T}\|_F}
\]
for each method to evaluate their performance. Meanwhile, we compare the average running time until algorithm  converges to some preset threshold.


\subsubsection{Data}
\begin{figure}
    \centering
	\begin{subfigure}[b]{0.32\linewidth}
	\includegraphics[width=\textwidth]
		{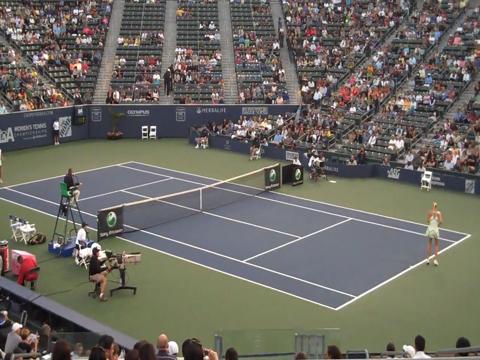}
		\caption{}
		\label{FIG:video1}
	\end{subfigure}
	\begin{subfigure}[b]{0.32\linewidth}
	\includegraphics[width=\textwidth]
		{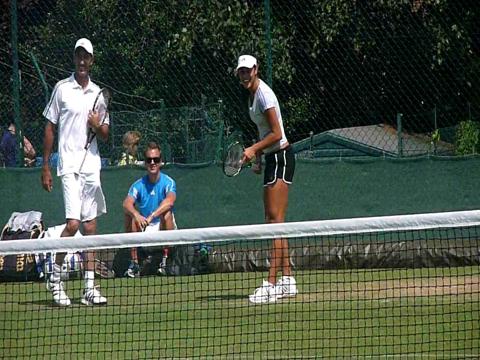}
		\caption{}
	\label{FIG:video2}
	\end{subfigure}
	\begin{subfigure}[b]{0.32\linewidth}
	\includegraphics[width=\textwidth]
		{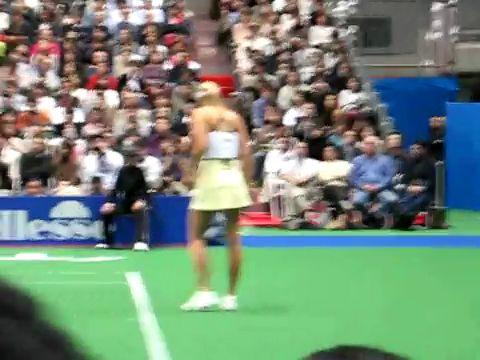}
		\caption{}
		\label{FIG:video3}
	\end{subfigure}
	  \caption{The first frame of tested videos}
	  \label{FIG: video}
\end{figure}
In this part, we have tested our algorithm on three video data, see \cite{FHL2018}. The video data are tennis-serve data from an Olympic Sports Dataset. The three videos are color videos. In our simulation, we use the same set-up like the one in \cite{FHL2018}, we pick 30 frames evenly from each video. For each image frame, the size is scaled to $360\times 480\times 3$. So each video is transformed into a 4-D tensor data of size $360\times 480\times 3\times 30$. The first frame of each video after preprocessed is illustrated in Figure \ref{FIG: video}.

\subsubsection{Numerical Results}
\begin{table}[!ht]
    \centering
    \caption{The Relative Square Error (RSE) and time spend for different algorithms on video data}
    \label{tab:video}
    \begin{tabular}{|c|c|c|c|c|}\hline
Video & Method & SR & RSE & Time(S) \\ \hline

\multirow{8}{*}{\includegraphics[width=0.12\textwidth]{img/1.jpg}}& wHOSVD-TV  & 10\% &  0.2080 &  82.26 \\
 &  &  30\% & 0.1418 & 50.40 \\
 &  &  50\% & 0.1045 & 41.31 \\
 &  &  80\% & 0.0566 & 33.21 \\
 \cline{2-5}& ST-HOSVD & 10\% &  N/A &  N/A \\
 &  &  30\% & 0.1941 & 521.94 \\
 &  &  50\% & 0.1381 & 175.82 \\
 &  &  80\% & 0.0667 & 128.68 \\\hline

\multirow{8}{*}{\includegraphics[width=0.12\textwidth]{img/2.jpg}} & wHOSVD-TV  & 10\% &  0.2694 &  35.21 \\
 &  &  30\% & 0.1888 & 21.08 \\
 &  &  50\% & 0.1411 & 16.55 \\
 &  &  80\% & 0.0767 & 12.88 \\
 \cline{2-5}
 & ST-HOSVD & 10\% &  N/A &  N/A \\
 &  &  30\% & 0.2249 & 1130.77 \\
 &  &  50\% & 0.1480 & 1304.31 \\
 &  &  80\% & 0.0749 & 976.65  \\\hline

\multirow{8}{*}{\includegraphics[width=0.12\textwidth]{img/3.jpg}} & wHOSVD-TV & 10\% &  0.2198 &  156.5 \\
 &  &  30\% & 0.1394 & 87.89 \\
 &  &  50\% & 0.0955 & 72.15 \\
 &  &  80\% & 0.0470 & 18.44 \\
 \cline{2-5}
 & ST-HOSVD & 10\% &  N/A &  N/A \\
 &  &  30\% & 0.1734 & 560.97 \\
 &  &  50\% & 0.1105 & 158.74 \\
 &  &  80\% & 0.0594 & 52.09 \\\hline
    \end{tabular}
\end{table}

\begin{figure}[h]
    \centering
    \includegraphics[width=0.5\textwidth]{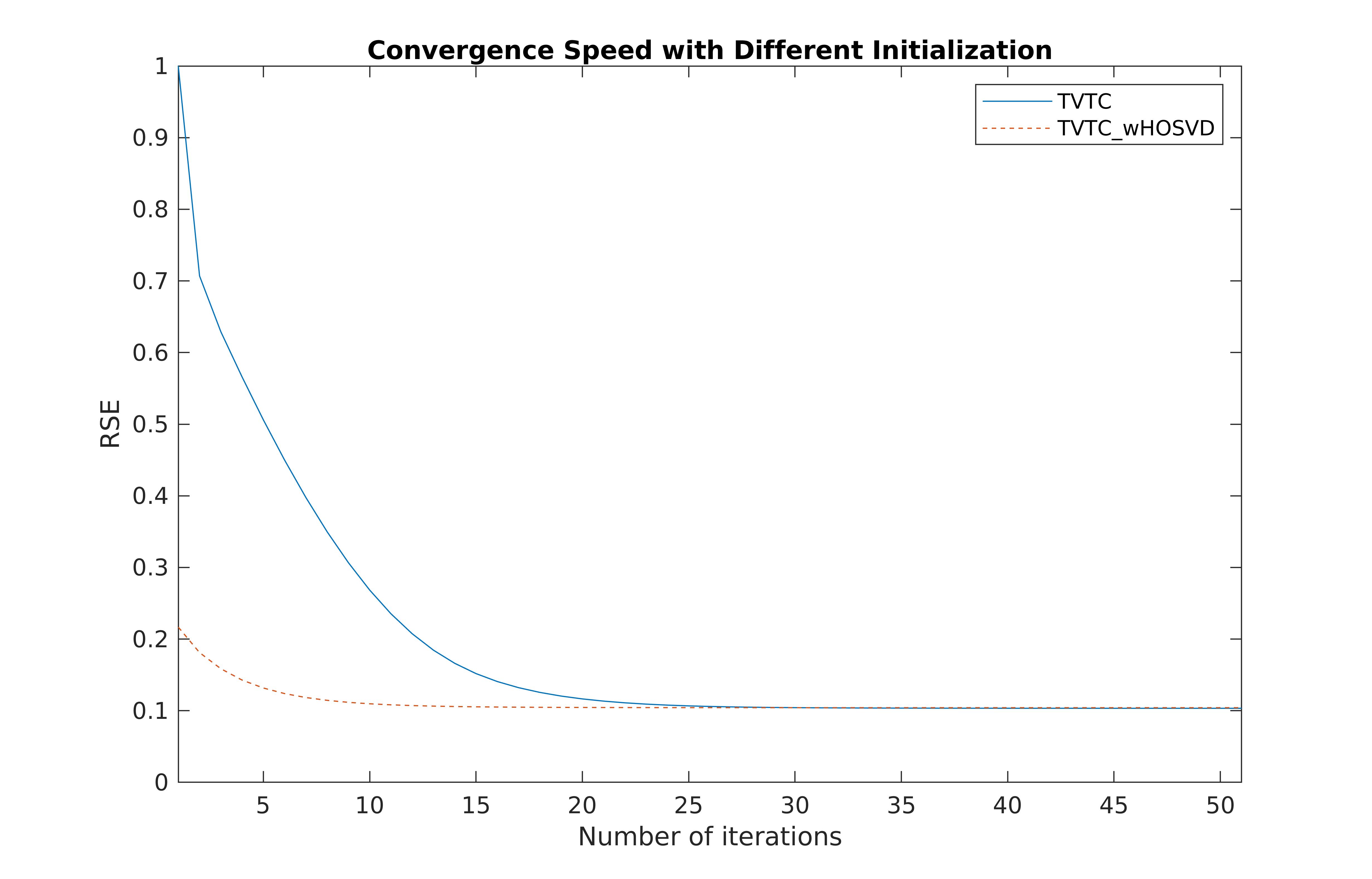}
     \caption{Comparison Between TVTC and wHOSVD-TV on video 1 with SR $=50\%$.}\label{fig:converge}
\end{figure}
The simulation results on the videos are reported in Table \ref{tab:video}  and Figure \ref{fig:converge}. In existent studies,  there are studies performed the same completion task on the same dataset (see \cite{FHL2018}.) In \cite{FHL2018}, the ST-HOSVD \cite{VVM2012}  had the best performance among several  low-rank based tensor completion algorithm.

We record the completion error and running time for each completion task and compare them with the previous low-rank based algorithms. One can observe that the TV-based algorithm is more compatible with video data most of the time. 

On the other hand, we have implemented total variations with zero filling initialization for the entries which are not observed and with the tensor obtained from weighted HOSVD which are termed TVTC and wHOSVD-TV respectively. The iterative results are shown in Figure \ref{fig:converge}, which shows that using the result from weighted HOSVD as initialization could notably reduce the iterations of TV-minimization for achieving the convergence threshold ($\|\tens{X}^k-\tens{X}^{k-1}\|_F<10^{-4}$).

\subsection{Discussion}

\begin{figure}[ht]
    \centering
    \begin{tikzpicture}[
 image/.style = {text width=0.5\linewidth, 
                 inner sep=0pt, outer sep=0pt},
node distance = 1mm and 1mm
                        ] 
\node [image] (frame1)
    {\includegraphics[width=\linewidth]{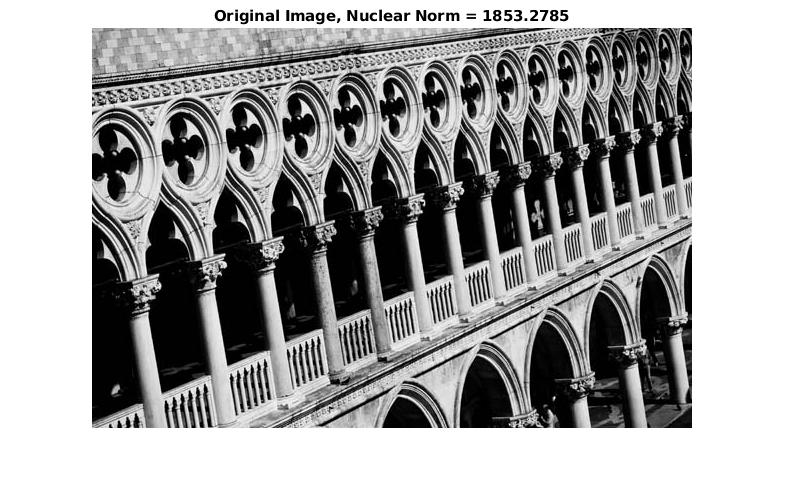}};
\node [image,right=of frame1] (frame2) 
    {\includegraphics[width=\linewidth]{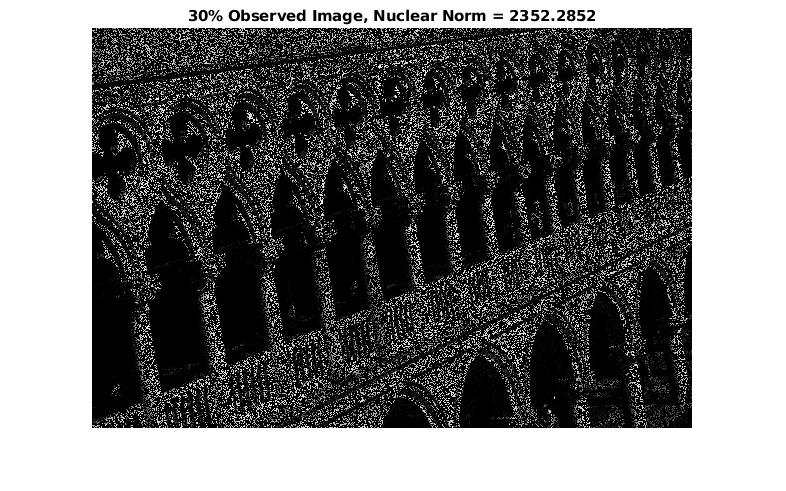}};
\node[image,below=of frame1] (frame3)
    {\includegraphics[width=\linewidth]{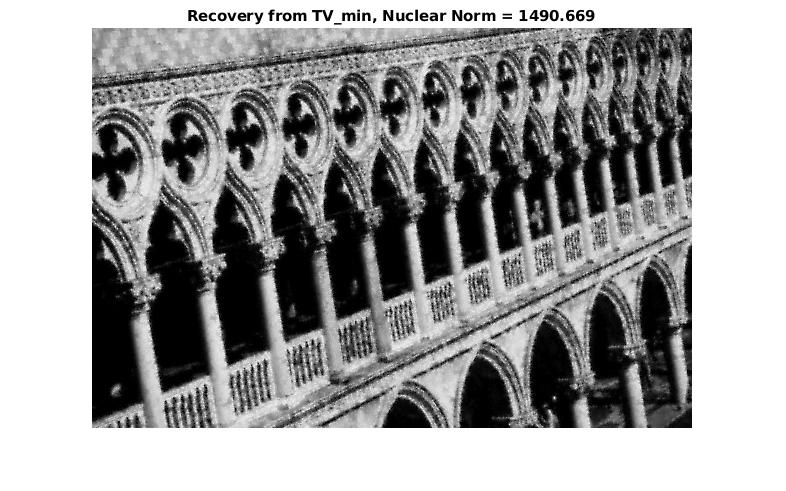}};
\node[image,right=of frame3] (frame4)
    {\includegraphics[width=\linewidth]{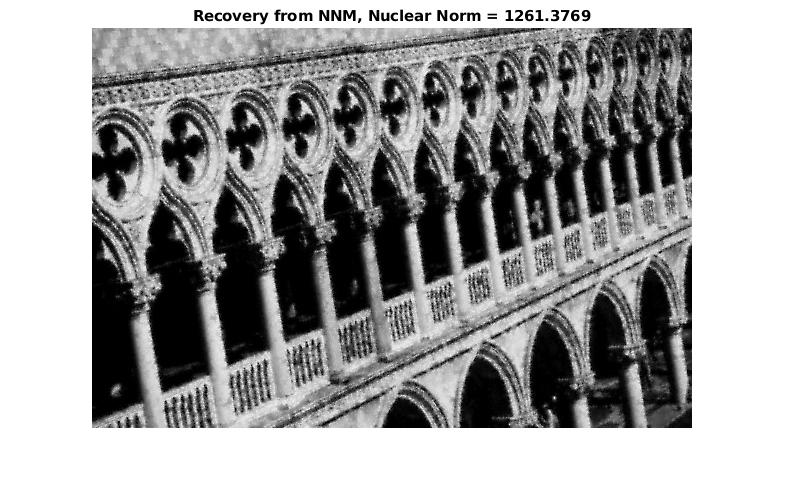}};
\node[image,below=of frame3] (frame5)
    {\includegraphics[width=\linewidth]{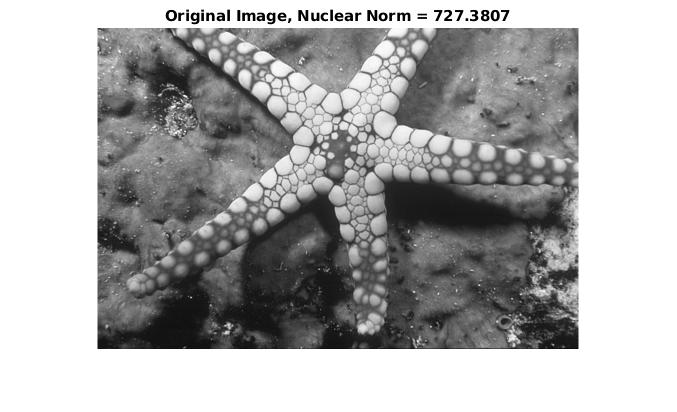}};
\node[image,right=of frame5] (frame6)
    {\includegraphics[width=\linewidth]{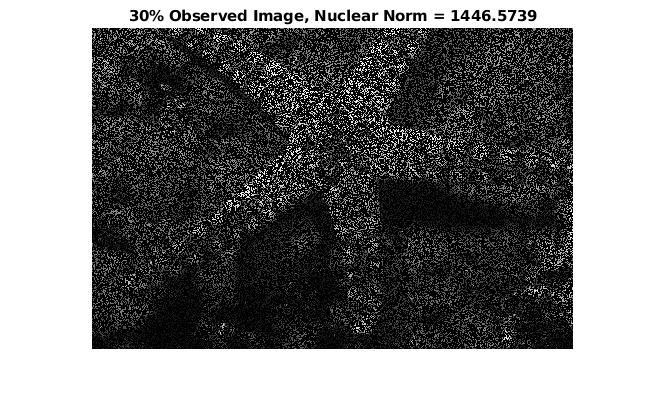}};
\node[image,below=of frame5] (frame7)
    {\includegraphics[width=\linewidth]{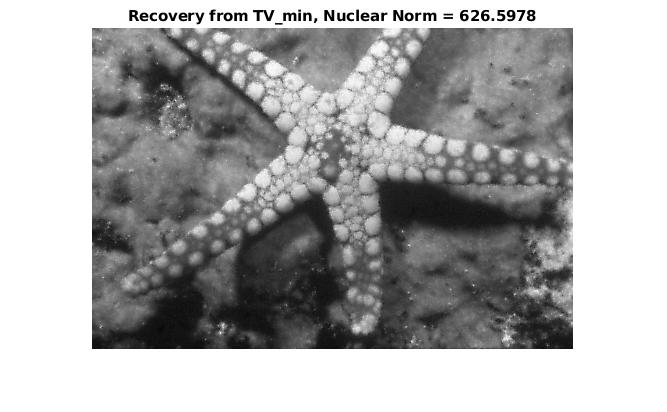}};
\node[image,right=of frame7] (frame8)
    {\includegraphics[width=\linewidth]{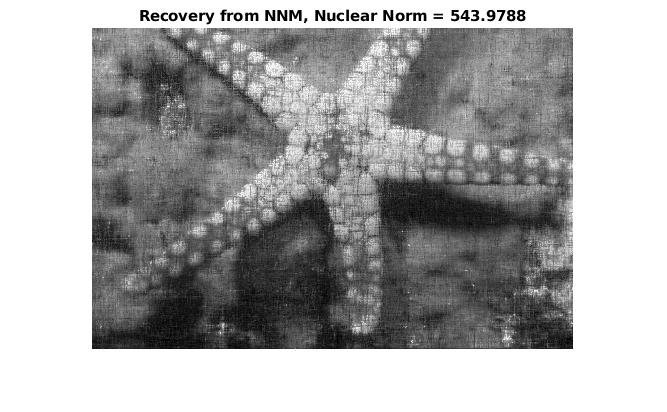}};
\end{tikzpicture}
    \caption{Nuclear Norm Comparison for Different Recovery Patterns}
    \label{FIG:NN}
\end{figure}
The relation between smoothness pattern and low-rank pattern is mysterious. When studying image-related data, both patterns are usually taken at the same time and converted to a mixed optimization structure. Through experiments, we find  that, with uniform random missing entries, result from total variation minimization on image-like data has singular values closer to the original image. 

We randomly mask 70\% entries for several grey-scale images and performed both nuclear norm minimization and total variation minimization on the images. Then the nuclear norms for original image, masked image, TV estimates, and NNM estimates are computed, see Figure \ref{FIG:NN}.
From this figure, we can see that the image recovered from TV minimization already gives a smaller nuclear norm compared to the original image, while NNM will bring this further away. By observing the singular values of the TV-recovered matrix and of the NNM recovered matrix, we can see that TV-minimization could better capture the smaller singular values, hence better preserve the overall structures of the original matrix.
\begin{figure}[h]
    \centering
    \includegraphics[width=0.4\textwidth]{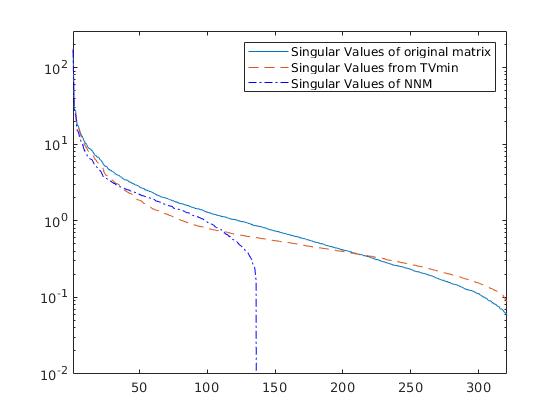}
    \caption{Comparison of singular values between TV-recovery and original image, original image is the same as in Figure \ref{FIG:NN}}
    \label{FIG:SV}
\end{figure}
Since both $\|\cdot\|_*$ and $\|\cdot\|_{TV}$ are convex functions, the mixed minimization problem with restricted observation entries
$$\widehat{\tens{X}}_{\text{mix}} = \min_\tens{X}\sum_i\alpha_i\|\tens{X}_{(i)}\|_* + \lambda \|\tens{X}\|_{TV}$$
will produces a result whose nuclear norm is between $\|\tens{X}_{NNM}\|_*$ and $\|\tens{X}_{TV}\|_*$, where $\widehat{\tens{X}}_{\text{NNM}}$ and $\widehat{\tens{X}}_{\text{TV}}$ are the results from each individual minimization problem (with the same constraint):
$$\widehat{\tens{X}}_{\text{NNM}} = \min_\tens{X}\sum_i\alpha_i\|\tens{X}_{(i)}\|_*,$$
$$\widehat{\tens{X}}_{\text{TV}} = \min_\tens{X}\|\tens{X}\|_{TV}.$$
Unlike user-rating data and synthetic low-rank tensor, the image-like data tends to have a non-trivial tail of singular values. Figure \ref{FIG:SV} shows the similarity between the original image and TV-recovered image, which gives us hints about the performance comparison between TV-recovery and low-rank recovery.
\section*{Acknowledgement}
Needell was partially supported by NSF CAREER \#1348721 and NSF BIGDATA \#1740325.

\bibliographystyle{plain} 
\bibliography{reference}

\end{document}